\documentclass[a4paper]{article}

\usepackage{hyperref}
\hypersetup{colorlinks, linkcolor=black, citecolor=black, urlcolor=black}
\usepackage{tikz}
\usepackage{url}
\usepackage{xcolor}

\usepackage{amsmath}
\usepackage{amssymb}
\usepackage{amsthm}

\newcommand{\conv}{\operatorname{conv}}
\newcommand{\R}{\mathbb{R}}
\newcommand{\xc}{\operatorname{xc}}
\newcommand{\zerovec}{\mathbb{O}}

\newtheorem{thm}{Theorem}

\newcounter{claimcounter}
\newcommand{\claim}[1]{%
    \stepcounter{claimcounter}
    \medskip
    \noindent
    \textit{Claim \theclaimcounter: #1}
    \medskip
}

\newcommand{\calR}{\mathcal{R}}

\colorlet{myblue}{blue!80!black}
\colorlet{mygreen}{green!50!black}
\colorlet{myred}{red!80!black}

\colorlet{blueTable}{myblue!30}
\colorlet{greenTable}{mygreen!30}
\colorlet{redTable}{red!30}

\newcommand{\blue}{{\color{myblue} blue} {}}
\newcommand{\green}{{\color{mygreen} green} {}}
\newcommand{\red}{{\color{myred} red} {}}
\newcommand{\blueR}{{\color{myblue} \mathcal{R}''}}
\newcommand{\redR}{{\color{myred} \mathcal{R}'}}

\newcommand{\femail}[1]{\href{mailto:#1}{\nolinkurl{#1}}}

\title{Extension complexities of Cartesian products involving a pyramid}
\author{Hans Raj Tiwary%
\thanks{KAM/ITI Charles University in Prague;
\femail{hansraj@kam.mff.cuni.cz};
Partially supported by project GA15-11559S of GA \v{C}R.}
\and Stefan Weltge%
\thanks{ETH Zurich;
\femail{weltges@ethz.ch}.}
\and Rico Zenklusen%
\thanks{ETH Zurich;
\femail{ricoz@math.ethz.ch};
Supported by the Swiss National Science Foundation grant 200021\_165866, ``New Approaches to Constrained Submodular
Maximization''.}
}
\date{}

\begin{document}
\maketitle
\begin{abstract}
    It is an open question whether the linear extension complexity of the Cartesian product of two polytopes $ P,Q $ is
    the sum of the extension complexities of $ P $ and $ Q $.
    We give an affirmative answer to this question for the case that one of the two polytopes is a pyramid.
\end{abstract}

\section{Introduction}
For a non-empty polytope $ P $, the \emph{linear extension complexity} of $ P $ is defined as the smallest number of
facets of any polytope that can be affinely projected onto $ P $, and is denoted by $ \xc(P) $.
Given any non-empty polytopes $P$ and $Q$, one can easily observe that $ \xc(P \times Q) \le \xc(P) + \xc(Q) $, while it is an open question whether this inequality actually holds as an equality, i.e., whether
\begin{equation}
    \label{eqProductEqualsSum}
    \xc(P \times Q) = \xc(P) + \xc(Q)
\end{equation}
holds in general.
This question has been asked at several occasions (see, e.g.,~\cite[Conj.~1]{GPS2016} or~\cite[Prob.~3]{W2016}) but it
seems that the most general case in which is it known that~\eqref{eqProductEqualsSum} holds is when one of the two
polytopes is a simplex.
The latter fact has been observed by several authors and can be explicitly found in~\cite[Cor.~10]{GPS2016}.
In this note, we prove that~\eqref{eqProductEqualsSum} holds whenever one of the two polytopes is a pyramid
(in Section~\ref{sec:prel} we recall the definition of a pyramid):
\begin{thm}
    \label{thmMain}
    Let $ P,Q $ be non-empty polytopes such that one of the two polytopes is a pyramid.
    Then we have $ \xc(P \times Q) = \xc(P) + \xc(Q) $.
\end{thm}
While pyramids are still very special polytopes, with respect to linear extensions they are closely related to their
bases, which can be arbitrary polytopes.
Indeed, given a pyramid $ P $ with base $ B $ it is easy to see that $ \xc(P) = \xc(B) + 1 $ holds.
Thus, although our proof crucially exploits the structure of Cartesian products involving a pyramid, we hope that our
result opens doors for further generalizations.

In the next section, we discuss basic ingredients needed for the proof of Theorem~\ref{thmMain} while the proof itself
is given in Section~\ref{secProof}.

\section{Preliminaries}\label{sec:prel}
A polytope $ P \subseteq \R^d $ is called a \emph{pyramid with base $ B \subseteq \R^d $ and apex $ v \in \R^d $} if $ P
= \conv(B \cup \{v\}) $ and $ v $ is not contained in the affine hull of $ B $.
Note that $ v $ is contained in every facet of $ P $ except for one which contains all remaining vertices of $ P $.

Let $ P = \{ x \in \R^d : \langle a_i,x \rangle \le b_i \ i=1,\dotsc,m \} = \conv \{ v_1,\dotsc,v_n \} $ for some $
a_1,\dotsc,a_m \in \R^d $, $ b_1,\dotsc,b_m \in \R $, and $ v_1,\dotsc,v_m \in \R^d $, where $ \langle \cdot,\cdot
\rangle $ denotes the Euclidean scalar product of $ \R^d $.
Then the matrix $ S \in \R^{m \times n}_{\ge 0} $ defined via $ S_{i,j} := b_i - \langle a_i,v_j \rangle $ is called a
\emph{slack matrix} of $ P $.
A well-known result of Yannakakis~\cite{Y1991} states that the linear extension complexity of $ P $ is equal to the
\emph{nonnegative rank} of $ S $, which is defined as the smallest number $ r_+(S) $ such that $ S $ can be written as
the sum of $ r_+(S) $ nonnegative rank-one matrices. The nonnegative rank $r_+(S)$ of a polytope is indeed well defined despite the fact its definition relies on the slack matrix $S$ which, in turn, is defined by a particular linear description of $P$. This follows from the fact that $r_+(S)$ neither depends on the scaling of the constraints used to describe $P$ nor on the potential presence of redundant constraints.

Although not needed for this work, the interested reader may consider the surveys~\cite{K2011,CCZ2013} and the book
chapter~\cite[Chap.~4]{CCZ2014} as excellent sources for background information and recent developments on linear
extended formulations.

In our proof, we make use of two simple facts about decompositions into nonnegative rank-one matrices:
Let $ S = R^1 + \dotsb + R^k $ where $ R^1,\dotsc,R^k $ are nonnegative rank-one matrices and suppose that $ S_{i,j} = 0
$ holds.
First, since all $ R^\ell $ are nonnegative, this implies $ (R^\ell)_{i,j} = 0 $ for all $ \ell $.
Second, since all $ R^\ell $ have rank one, for every pair of indices $ (i',j') $ and every $ \ell $ we must have $
(R^\ell)_{i',j} = 0 $ or $ (R^\ell)_{i,j'} = 0 $.

Given two polytopes $ P,Q $ with
\begin{align*}
    P & = \{ x \in \R^{d_P} : \langle a^P_i,x \rangle \le b^P_i \ i=1,\dotsc,m_P \} = \conv \{ v^P_1,\dotsc,v^P_{n_P} \}
    \intertext{and}
    Q & = \{ y \in \R^{d_Q} : \langle a^Q_i,y \rangle \le b^Q_i \ i=1,\dotsc,m_Q \} = \conv \{ v^Q_1,\dotsc,v^Q_{n_Q} \},
\end{align*}
one immediately obtains
\begin{alignat*}{10}
    P \times Q & = \{ (x,y) \in \R^{d_P} \times \R^{d_Q} \ \, : & & \langle a^P_i,x \rangle \le b^P_i \ i=1,\dotsc,m_P, \\
               &                                           & & \langle a^Q_i,y \rangle \le b^Q_i \ i=1,\dotsc,m_Q \} \\
               & = \conv \{ (v_i^P, v_j^Q) : i \in [n_P], & & \, j \in [n_Q] \}.
\end{alignat*}
Thus, if $ S \in \R^{m_P \times n_P}_{\ge 0} $ and $ T = [ t_1 \dotsb t_{n_Q} ] \in \R^{m_Q \times n_Q}_{\ge 0}
$ are slack matrices of $ P $ and $ Q $, respectively, then the matrix
\medskip
\begin{center}
    \begin{tikzpicture}[scale=2]
        \draw (0.5,1) -- (4.5,1);
        \draw (0.5,1.5) -- (4.5,1.5);
        \draw (0.5,2.5) -- (4.5,2.5);

        \draw (0.5,1) -- (0.5,2.5);
        \draw (1.5,1) -- (1.5,2.5);
        \draw (2.5,1) -- (2.5,2.5);
        \draw (3.5,1) -- (3.5,2.5);
        \draw (4.5,1) -- (4.5,2.5);

        \node at (1,2) {$ S $};
        \node at (2,2) {$ S $};
        \node at (3,2) {$ \dotsb $};
        \node at (4,2) {$ S $};

        \node at (1,1.25) {$ t_1 \dotsb t_1 $};
        \node at (2,1.25) {$ t_2 \dotsb t_2 $};
        \node at (3,1.25) {$ \dotsb $};
        \node at (4,1.25) {$ t_{n_Q} \dotsb t_{n_Q} $};

        \node at (5.5,1.75) {$ \in \R^{(m_P + m_Q) \times (n_P \cdot n_Q)}_{\ge 0} $};
    \end{tikzpicture}
\end{center}
\medskip
is a slack matrix of $ P \times Q $, where $ t_1,\dotsc,t_{n_Q} \in \R^{m_Q}_{\ge 0} $ denote the columns of $ T $.
The columns of the above slack matrix correspond, from left to right, to the vertices $(v_1^P, v_1^Q), (v_2^P,v_1^Q),\ldots, (v_{n_p}^P,v_1^Q), (v_1^P,v_2^Q), (v_2^P,v_2^Q), \ldots ,(v_{n_p}^P,v_{n_Q}^Q)$. Moreover, the first block of rows correspond to the constraints of $P$ and the second block of rows to the constraints of $Q$.

\section{Proof of Theorem~\ref{thmMain}}
\label{secProof}
We may assume that $ Q $ is a pyramid.
First, note that there exists a slack matrix $ S \in \R^{m_P \times n_P}_{\ge 0} $ of $ P $ such that every row contains
at least one entry being zero.
Indeed, every row containing no entry being zero corresponds to a redundant inequality and hence can be removed from the
description of $ P $.
Second, by assuming that the description of $Q$ does not contain any redundant inequalities, the slack matrix $T\in \R^{m_Q\times n_Q}$ of $Q$ has the form
%
%
\medskip
\begin{center}
    \begin{tikzpicture}[scale=1]
        \draw (0.5,1) -- (2,1);
        \draw (0.5,1.5) -- (2,1.5);
        \draw (0.5,2.5) -- (2,2.5);

        \draw (0.5,1) -- (0.5,2.5);
        \draw (1.5,1) -- (1.5,2.5);
        \draw (2,1) -- (2,2.5);

        \node at (0,1.75) {$ T = $};

        \node at (1,2) {$ T' $};
        \node at (1.75,2) {$ \zerovec $};

        \node at (1,1.25) {$ \zerovec $};
        \node at (1.75,1.25) {$ 1 $};
    \end{tikzpicture}
\end{center}
where $ T' \in \R^{(m_Q - 1) \times (n_Q - 1)}_{\ge 0} $.
Thus, the matrix $ A \in \R^{(m_P + m_Q) \times (n_P \cdot n_Q)}_{\ge 0} $ defined via
\medskip
\begin{center}
    \begin{tikzpicture}[scale=2]
        \fill[greenTable] (0.5,1.5) rectangle (4.5,2.5);
        \fill[blueTable] (4.5,1.5) rectangle (5.5,2.5);
        \fill[redTable] (0.5,1) rectangle (4.5,1.5);
        \fill[blueTable] (4.5,0.7) rectangle (5.5,1);

        \draw (0.5,0.70) -- (5.5,0.70);
        \draw (0.5,1) -- (5.5,1);
        \draw (0.5,1.5) -- (5.5,1.5);
        \draw (0.5,2.5) -- (5.5,2.5);

        \draw (0.5,0.7) -- (0.5,2.5);
        \draw (1.5,0.7) -- (1.5,2.5);
        \draw (2.5,0.7) -- (2.5,2.5);
        \draw (3.5,0.7) -- (3.5,2.5);
        \draw (4.5,0.7) -- (4.5,2.5);
        \draw (5.5,0.7) -- (5.5,2.5);

        \node at (0.2,1.5) {$ A := $};

        \node at (1,2) {$ S $};
        \node at (2,2) {$ S $};
        \node at (3,2) {$ \dotsb $};
        \node at (4,2) {$ S $};
        \node at (5,2) {$ S $};

        \node at (1,1.25) {$ t'_1 \dotsb t'_1 $};
        \node at (2,1.25) {$ t'_2 \dotsb t'_2 $};
        \node at (3,1.25) {$ \dotsb $};
        \node at (4,1.25) {$ t'_k \dotsb t'_k $};
        \node at (5,1.25) {$ \zerovec $};

        \node at (1,0.85) {$ \zerovec $};
        \node at (2,0.85) {$ \zerovec $};
        \node at (3,0.85) {$ \dotsb $};
        \node at (4,0.85) {$ \zerovec $};
        \node at (5,0.85) {$ 1 \dotsb 1 $};
    \end{tikzpicture}
\end{center}
\medskip
is a slack matrix of $ P \times Q $, where $ t'_1,\dotsc,t'_k \in \R^{m_Q - 1}_{\ge 0} $ are the columns of $ T'
$ (here $ k = n_Q - 1 $).
Recall that we have $ \xc(P \times Q) = r_+(A) $, $ \xc(P) = r_+(S) $, and $ \xc(Q) = r_+(T) $.
Furthermore, it is straightforward to check that $ r_+(T) = r_+(T') + 1 $ holds.
Thus, it remains to show that
\[
    r_+(A) \ge r_+(S) + r_+(T') + 1
\]
holds.
For the sake of contradiction, let us assume that we have
\[
    r_+(A) \le r_+(S) + r_+(T'),
\]
i.e., there exists a set $ \calR $ of nonnegative rank-one matrices in $ \R^{(m_P + m_Q) \times (n_P \cdot n_Q)}_{\ge 0}
$ with $ |\calR| \le r_+(S) + r_+(T') $ whose sum is equal to $ A $.
Let $ \redR $ and $ \blueR $ denote the set of matrices in $ \calR $ that have support in the \red and \blue
parts of $ A $, respectively.

\claim{The sets $ \redR $ and $ \blueR $ form a partition of $ \calR $ satisfying $ |\redR| = r_+(T') $ and $ |\blueR| =
r_+(S) $.}

First, observe that $ \redR $ and $ \blueR $ are disjoint due to the $\zerovec$-block within $A$ that is below the \blue $S$-block.
Since the \red part of $ A $ contains $ T' $ as a submatrix, we must have $ |\redR| \ge r_+(T') $, and since the \blue
part contains $ S $ as a submatrix, we must have $ |\blueR| \ge r_+(S) $, which yields the claim.

\claim{There exists at least one matrix in $ \redR $ that has support in the \green part of $ A $.}

Since the nonnegative rank of the \green submatrix of $ A $ is equal to the nonnegative rank of $ S $, at least $ r_+(S)
$ matrices in $ \calR $ must have support in this part.
Note that at least one matrix in $ \blueR $ has support in the last row of the \blue part of $ A $ and hence it cannot
have support in the \green part of $ A $.
The claim follows since $ |\blueR| = r_+(S) $.

\claim{Let $ R \in \redR $ and pick exactly one column of each of the $ k $ \red submatrices of $ A $. Then $ R $ has
support in at least one of these columns.}

Suppose the contrary.
Then we can pick exactly one column of each of the $ k $ \red submatrices of $ A $ such that $ R $ has no
support on any of these columns.
Restricting to the submatrix formed by these columns, observe that this submatrix is identical to $ T' $ but can be
written as the sum of all matrices in $ \redR \setminus \{R\} $ and hence $ r_+(T') \le |\redR| - 1 = r_+(T') - 1 $, a
contradiction.

\claim{No matrix in $ \redR $ can have support in the \green part of $ A $ (a contradiction to Claim~2).}

Assume that there is some $ R \in \redR $ that has a positive entry $ e_1 $ in the \green part of $ A $.
By our choice of $ S $, every of the first $ k $ blocks of $ A $ contains a column of $ A $ in which this row has a zero
entry.
By the previous claim, $ R $ has a positive entry $ e_2 $ in the \red part of one of these columns.
Restricting $ R $ to the two-by-two submatrix containing the entries $ e_1,e_2 $, it looks as follows (up to swapping
its columns):
\medskip
\begin{center}
    \begin{tikzpicture}
        \fill[greenTable] (0.5,1.5) rectangle (2.5,2.5);
        \fill[redTable] (0.5,0.5) rectangle (2.5,1.5);

        \node at (1,2) {\footnotesize $ e_1 > 0 $};
        \node at (2,2) {\footnotesize $ 0 $};
        \node at (2,1) {\footnotesize $ e_2 > 0 $};
        \node at (1,1) {\footnotesize $ * $};

        \draw (0.5,0.5) -- (2.5,0.5);
        \draw (0.5,1.5) -- (2.5,1.5);
        \draw (0.5,2.5) -- (2.5,2.5);
        \draw (0.5,0.5) -- (0.5,2.5);
        \draw (1.5,0.5) -- (1.5,2.5);
        \draw (2.5,0.5) -- (2.5,2.5);
    \end{tikzpicture}
\end{center}
\medskip
However, there is no rank-one matrix with such a sign pattern. \qed

\bibliographystyle{plain}
\bibliography{references}

\end{document}